\documentclass[11pt]{amsart}

\usepackage[all]{xy}
\usepackage{color, pstricks}
\usepackage{mathrsfs}
\usepackage{amsmath,amsthm, amssymb, amsgen,amsxtra,amsfonts, amsbsy} 
\usepackage[all]{xy}
\usepackage{verbatim}

\usepackage{times}
\newcommand{\mathds}[1]{{\mathbb #1}}

\begin{document}
%
%   D e f i n i t i o n s
%
%
\theoremstyle{definition}
\newtheorem{Definition}{Definition}[section]
\newtheorem*{Definitionx}{Definition}
\newtheorem{Convention}{Definition}[section]
\newtheorem{Construction}{Construction}[section]
\newtheorem{Example}[Definition]{Example}
\newtheorem{Examples}[Definition]{Examples}
\newtheorem{Remark}[Definition]{Remark}
\newtheorem*{Remarkx}{Remark}
\newtheorem{Remarks}[Definition]{Remarks}
\newtheorem{Caution}[Definition]{Caution}
\newtheorem{Conjecture}[Definition]{Conjecture}
\newtheorem*{Conjecturex}{Conjecture}
\newtheorem{Question}[Definition]{Question}
\newtheorem{Questions}[Definition]{Questions}
\newtheorem*{Acknowledgements}{Acknowledgements}
\newtheorem*{Organization}{Organization}
\newtheorem*{Disclaimer}{Disclaimer}
\theoremstyle{plain}
\newtheorem{Theorem}[Definition]{Theorem}
\newtheorem*{Theoremx}{Theorem}
\newtheorem{Proposition}[Definition]{Proposition}
\newtheorem*{Propositionx}{Proposition}
\newtheorem{Lemma}[Definition]{Lemma}
\newtheorem{Corollary}[Definition]{Corollary}
\newtheorem*{Corollaryx}{Corollary}
\newtheorem{Fact}[Definition]{Fact}
\newtheorem{Facts}[Definition]{Facts}
\newtheoremstyle{voiditstyle}{3pt}{3pt}{\itshape}{\parindent}%
{\bfseries}{.}{ }{\thmnote{#3}}%
\theoremstyle{voiditstyle}
\newtheorem*{VoidItalic}{}
\newtheoremstyle{voidromstyle}{3pt}{3pt}{\rm}{\parindent}%
{\bfseries}{.}{ }{\thmnote{#3}}%
\theoremstyle{voidromstyle}
\newtheorem*{VoidRoman}{}

% abgeschrieben aus The LaTeX Companion, 2nd edition,
% von Mittelback & Goossens
%
\newcommand{\prf}{\par\noindent{\sc Proof.}\quad}
\newcommand{\blowup}{\rule[-3mm]{0mm}{0mm}}
\newcommand{\cal}{\mathcal}
\newcommand{\Aff}{{\mathds{A}}}
\newcommand{\BB}{{\mathds{B}}}
\newcommand{\CC}{{\mathds{C}}}
\newcommand{\EE}{{\mathds{E}}}
\newcommand{\FF}{{\mathds{F}}}
\newcommand{\GG}{{\mathds{G}}}
\newcommand{\HH}{{\mathds{H}}}
\newcommand{\NN}{{\mathds{N}}}
\newcommand{\ZZ}{{\mathds{Z}}}
\newcommand{\PP}{{\mathds{P}}}
\newcommand{\QQ}{{\mathds{Q}}}
\newcommand{\RR}{{\mathds{R}}}
\newcommand{\Liea}{{\mathfrak a}}
\newcommand{\Lieb}{{\mathfrak b}}
\newcommand{\Lieg}{{\mathfrak g}}
\newcommand{\Liem}{{\mathfrak m}}
\newcommand{\ideala}{{\mathfrak a}}
\newcommand{\idealb}{{\mathfrak b}}
\newcommand{\idealg}{{\mathfrak g}}
\newcommand{\idealm}{{\mathfrak m}}
\newcommand{\idealp}{{\mathfrak p}}
\newcommand{\idealq}{{\mathfrak q}}
\newcommand{\idealI}{{\cal I}}
\newcommand{\lin}{\sim}
\newcommand{\num}{\equiv}
\newcommand{\dual}{\ast}
\newcommand{\iso}{\cong}
\newcommand{\homeo}{\approx}
\newcommand{\mm}{{\mathfrak m}}
\newcommand{\pp}{{\mathfrak p}}
\newcommand{\qq}{{\mathfrak q}}
\newcommand{\rr}{{\mathfrak r}}
\newcommand{\pP}{{\mathfrak P}}
\newcommand{\qQ}{{\mathfrak Q}}
\newcommand{\rR}{{\mathfrak R}}
\newcommand{\bmu}{\boldsymbol{\mu}}
\newcommand{\balpha}{\boldsymbol{\alpha}}
%
%  evtl. auch \"uber \mathbb oder \Bbb
%
\newcommand{\OO}{{\cal O}}
\newcommand{\numero}{{n$^{\rm o}\:$}}
\newcommand{\mf}[1]{\mathfrak{#1}}
\newcommand{\mc}[1]{\mathcal{#1}}
\newcommand{\into}{{\hookrightarrow}}
\newcommand{\onto}{{\twoheadrightarrow}}
\newcommand{\Spec}{{\rm Spec}\:}
\newcommand{\BigSpec}{{\rm\bf Spec}\:}
\newcommand{\Spf}{{\rm Spf}\:}
\newcommand{\Proj}{{\rm Proj}\:}
\newcommand{\Pic}{{\rm Pic}}
\newcommand{\Br}{{\rm Br}}
\newcommand{\NS}{{\rm NS}}
\newcommand{\Sym}{{\mathfrak S}}
\newcommand{\Aut}{{\rm Aut}}
\newcommand{\Autp}{{\rm Aut}^p}
\newcommand{\Hom}{{\rm Hom}}
\newcommand{\Ext}{{\rm Ext}}
\newcommand{\ord}{{\rm ord}}
\newcommand{\coker}{{\rm coker}\,}
\newcommand{\divisor}{{\rm div}}
\newcommand{\Def}{{\rm Def}}
\newcommand{\piet}{{\pi_1^{\rm \acute{e}t}}}
\newcommand{\Het}[1]{{H_{\rm \acute{e}t}^{{#1}}}}
\newcommand{\Hfl}[1]{{H_{\rm fl}^{{#1}}}}
\newcommand{\Hcris}[1]{{H_{\rm cris}^{{#1}}}}
\newcommand{\HdR}[1]{{H_{\rm dR}^{{#1}}}}
\newcommand{\hdR}[1]{{h_{\rm dR}^{{#1}}}}
\newcommand{\defin}[1]{{\bf #1}}
\newcommand{\oX}{\cal{X}}
\newcommand{\oA}{\cal{A}}
\newcommand{\oY}{\cal{Y}}

\title[The Picard Rank of an Enriques Surface]{The Picard Rank of an Enriques Surface}
\author{Christian Liedtke}
\address{TU M\"unchen, Zentrum Mathematik - M11, Boltzmannstr. 3, D-85748 Garching bei M\"unchen, Germany}
\curraddr{}
\email{liedtke@ma.tum.de}

\date{\today}
\subjclass[2010]{14J28, 14C22, 14G17}

\begin{abstract}
   In this note, we use crystalline methods and the Tate-conjecture to give a
   short proof that the Picard rank of an Enriques surface 
   is equal to its second Betti number.
\end{abstract}

\maketitle

\section{Introduction}

Enriques surfaces are one of the four classes of minimal, smooth, and proper 
surfaces of Kodaira dimension zero.
The following fundamental result relates the Picard rank $\rho$
to the second Betti number $b_2$ of these surfaces.

\begin{Theorem}[Bombieri--Mumford \cite{bm3}]
  \label{b2=rho=10}
  Let $X$ be an Enriques surface over an algebraically closed field $k$.
  Then, $\rho(X)=b_2(X)=10$.
\end{Theorem}

Using this result, it is not difficult to show 
that the N\'eron--Severi lattice of an Enriques surface is even, 
unimodular, of signature $(1,9)$, and of discriminant $-1$, see \cite[Corollaire II.7.3.7]{Illusie}.
Thus, it is isometric to $U\perp E_8$ by lattice theory, see \cite[Chapter I.5]{CDL}.
In particular, there exist non-zero isotropic vectors, which implies
that every Enriques surface carries a genus-one fibration.
Moreover, this result is also essential for the analysis of 
linear systems \cite{CossecPicard}, projective models \cite{CossecModels}, \cite{Liedtke}, 
automorphism groups \cite{Barth}, and moduli spaces \cite{Gritsenko} of these surfaces.
\medskip

If $k=\CC$, then Theorem \ref{b2=rho=10} is an easy consequence
of $H^2(\OO_X)=0$ and the Lefschetz theorem on $(1,1)$ classes.
On the other hand, the known proofs of this result if
${\rm char}(k)>0$ are rather delicate and complicated.
\begin{enumerate}
 \item The first proof is due to Bombieri and Mumford \cite{bm3}, where they first establish 
  with some effort the existence of a genus-one fibration $f:X\to\PP^1$.
  Using this, they determine $\rho(X)$ via passing to the Jacobian surface $J(X)\to\PP^1$ 
  of $f$, which is a rational surface, and thus, satisfies $\rho=b_2$.
 \item Another proof is due to Lang \cite{Lang}, who first establishes lifting of $X$ 
   to characteristic zero for some classes of Enriques surfaces and then, 
   he uses the result in characteristic zero and specialization arguments.
   In the remaining cases, where lifting was unclear, he proves that $X$ is unirational,
   and then, uses results of Shioda to conclude.
\end{enumerate}
\medskip

In this note, we give a conceptual proof of Theorem \ref{b2=rho=10} 
that neither makes heavy use of special properties of Enriques surfaces,
nor relies on case-by-case analyses.
The idea of our proof is similar to the easy proof over the complex numbers:
we merely use that the Witt-vector cohomology group $H^2(W\OO_X)$ 
is torsion
(note that $H^2(\OO_X)$ may be non-zero in positive characteristic),
as well as the Tate-conjecture for Enriques surfaces over finite fields,
which is an arithmetic analog of the Lefschetz theorem on $(1,1)$ classes.
We refer to Remark \ref{rem: analog} for details.
\medskip

This note is organized as follows:

In Section \ref{sec: assuming Tate}, we give a short proof of 
Theorem \ref{b2=rho=10} assuming the Tate-conjecture 
for Enriques surfaces over finite fields.

In order to obtain an unconditional proof, we establish
in Section \ref{sec: Tate} the Tate-conjecture for Enriques surfaces
over finite fields, using as little special properties of these surfaces
as possible.

 \begin{Acknowledgements}
  I thank Igor Dolgachev for comments and discussion, as well as
  the referee for comments and careful proof-reading.
 \end{Acknowledgements}

\section{A short proof assuming the Tate-conjecture}
\label{sec: assuming Tate}

In this section, we first recall the definition of Enriques surfaces,
as well as a couple of their elementary properties.
Then, we reduce Theorem \ref{b2=rho=10} to the case of finite fields, and 
finally, give a short proof of Theorem \ref{b2=rho=10}
assuming the Tate-conjecture for Enriques surfaces over finite fields.

\subsection{Enriques surfaces}
Let $X$ be a smooth and proper variety (geometrically integral scheme
of finite type) over a field $k$.
We denote numerical equivalence of divisors on $X$ by $\equiv$
and define the $i$.th Betti number $b_i$ of $X$ to be the $\QQ_\ell$-dimension of 
$\Het{i}(X,\QQ_\ell)$, where $\ell$ is a prime different from ${\rm char}(k)$.
For a fixed algebraic closure $\overline{k}$ of $k$, we set
$\overline{X}:=X\times_{\Spec k}\Spec\overline{k}$.

\begin{Definition}
  A smooth and proper surface $X$ over an algebraically closed field $k$ 
  is called an \emph{Enriques surface} if
  $$
     \omega_X\num\OO_X \mbox{ \quad and \quad } b_2(X)\,=\,10.
  $$
  Moreover, if $k$ is an arbitary field,
  then a smooth and proper variety $X$ over $k$ is called an
  Enriques surface if $\overline{X}$ is an Enriques surface over $\overline{k}$.
\end{Definition}

From the table in the introduction of \cite{bm2}, we obtain the following
equalities and bounds on the cohomology of Enriques surfaces
\begin{equation}
 \label{invariants}
    b_1(X)\,=\,0, \mbox{ \quad } b_2(X)\,=\,10, \mbox{ \quad and \quad }h^1(\OO_X)\,=\,h^2(\OO_X)\,\leq\,1.
\end{equation}
This is actually everything needed to prove Theorem \ref{thm: main assuming Tate} below.
We remark that Enriques surfaces with $h^2(\OO_X)\neq0$ do exist 
in characteristic $2$, see \cite{bm3}.

\subsection{Slope one and reduction to the case of finite fields}
Let $W=W(k)$ be the Witt ring of a perfect field $k$ and let
$K$ be the field of fractions of $W$.
Let $X$ be a smooth and proper variety over $k$.
Then, $b_i(X)$ is equal to the rank of the
$W$-module $\Hcris{i}(X/W)$.
The following is a straight forward generalization of
\cite[Proposition II.7.3.2]{Illusie}.

\begin{Proposition}
 \label{prop: h2w is torsion}
  Let $X$ be a smooth and projective variety over an algebraically closed
  field $k$ of positive characteristic that satisfies
  $$
     \frac{1}{2}b_1(X)=h^1(X,\,\OO_X)-h^2(X,\,\OO_X).
 $$
 Then, the $F$-isocrystal $\Hcris{2}(X/W)\otimes_W K$ is of slope one and
 $$
   H^2(X,W\OO_X) \,=\, H^2(X,W\OO_X)_{\rm tors} \,=\,
   H^2(X,W\OO_X)_{V-{\rm tors}}\,,
 $$ 
 where ${\rm tors}$ denotes torsion as $W$-module and $V-{\rm tors}$
 denotes $V$-torsion.
\end{Proposition}

\prf
By \cite[Remarque II.6.4]{Illusie}, the $V$-torsion $H^2_{V-{\rm tors}}$
of $H^2(W\OO_X)$ is isomorphic to ${\rm DM}(\Pic^0_{X/k}/\Pic^0_{X/k,{\rm red}})$,
where ${\rm M}(-)$ denotes the contravariant Dieudonn\'e module and
${\rm D}(-)={\rm Hom}_W(-,K/W)$.
Thus, by Dieudonn\'e theory, the $k$-dimension of 
$H^2_{V-{\rm tors}}/VH^2_{V-{\rm tors}}$ is equal to the dimension 
of the Zariski tangent space of $\Pic^0_{X/k}/\Pic^0_{X/k,{\rm red}}$,
which is equal to $h^1(\OO_X)-\frac{1}{2}b_1(X)$.
Thus, in the exact sequence
$$
...\,\to\,H^1(\OO_X)\,\to\,H^2(W\OO_X)
\,\stackrel{V}{\longrightarrow}\,H^2(W\OO_X) 
\,\stackrel{\alpha}{\longrightarrow}\, H^2(\OO_X)\,\to\,...,
$$
the restriction $\alpha|_{H^2_{V-{\rm tors}}}:H^2_{V-{\rm tors}}\to H^2(\OO_X)$ is
surjective by our assumptions.
Next, we set $L:=H^2(W\OO_X)/H^2_{V-{\rm tors}}$ and 
denote the map induced by $V$ on $L$ again by $V$.
Using the snake lemma, we conclude $L/VL=0$.
As explained in the proof of \cite[Proposition II.7.3.2]{Illusie}, $L$
is $V$-adically separated, which implies $L=0$.
Thus, $H^2(W\OO_X)=H^2_{V-{\rm tors}}$ and
this $W$-module is torsion.

Since the slope spectral sequence of $X$ degenerates 
up to torsion \cite[Th\'eor\`eme II.3.2]{Illusie}, we conclude
$$
   0\,=\,H^2(W\OO_X)\otimes_{W}K \,=\, 
   \left( \Hcris{2}(X/W)\otimes_{W}K \right)_{[0,1[}\,.
$$   
Since $X$ is projective over $k$, the hard Lefschetz theorem 
(see \cite{IllusieReport} or the discussion in \cite[Section II.5.B]{Illusie})
implies that also the part of slope $]1,2]$ is zero.
Thus, $\Hcris{2}(X/W)\otimes_{W}K$ is of slope one.
\qed\medskip

\begin{Proposition}[Ekedahl--Hyland--Shepherd-Barron]
 \label{prop: EHSB}
 Let $f:{\cal X}\to S$ be a smooth and projective morphism such
 that $S$ is Noetherian, $f_\ast\OO_{\cal X}\iso\OO_S$, and such that
 $$
     \frac{1}{2}b_1({\cal X}_{\bar{s}})=h^1(\OO_{{\cal X}_{\bar{s}}})-h^2(\OO_{{\cal X}_{\bar{s}}})
 $$
 for every geometric point $\bar{s}\to S$.
 Then, the geometric Picard rank in this family is locally constant.
\end{Proposition}

\prf
This is a special case of \cite[Proposition 4.2]{EHSB}.
\qed\medskip

\begin{Corollary}
  \label{cor: reduce to finite}
  In order to prove Theorem \ref{b2=rho=10}, it suffices to establish it for
  Enriques surfaces that can be defined over finite fields.
\end{Corollary}

\prf
Let $X$ be an Enriques surface over an algebraically closed field $k$.
Then, there exists a sub-$\ZZ$-algebra $R$ of $k$ that is 
of finite type over $\ZZ$ and a smooth and projective morphism
${\cal X}\to S:=\Spec R$ with ${\cal X}\times_S\Spec k\iso X$.
Moreover, if $s\in S$ is a closed point, then the residue field 
$\kappa(s)$ is a finite field.
% (See \cite[Section 6]{IllusieFrobenius} for details.)
In particular, the geometric fiber ${\cal X}_{\bar{s}}$ 
is an Enriques surface over $\overline{\kappa(s)}$ and
we have $\rho({\cal X}_{\bar{s}})=b_2({\cal X}_{\bar{s}})=10$
by assumption.
Using Proposition \ref{prop: EHSB}, the assertion follows.
\qed\medskip

\subsection{The Tate-conjecture (for divisors over finite fields)}
\label{subsec: Tate conjecture}
Let $X$ be a smooth and proper variety of dimension $d$ 
over a finite field $\FF_q$, let $N_r(X)$ to be the number of
$\FF_{q^r}$-rational points of $X$, and let
$$
 Z(X,t) \,:=\, \exp \left( \sum_{r=1}^\infty N_r(X)\, \frac{t^r}{r} \right)
 \,=\,
 \frac{P_1(t)\cdot P_3(t)\cdots P_{2d-1}(t)}{P_0(t)\cdot P_2(t)\cdots P_{2d}(t)}
$$
be the zeta function of $X$ as in \cite{Deligne}.
By loc.cit., there exist $\alpha_i\in\overline{\QQ}$ such that
\begin{equation}
 \label{P2}
  P_2(t) \,=\, \prod_{i=1}^{b_2(X)} (1-\alpha_i t)
\end{equation}
and such that for every embedding of fields $\QQ(\alpha_i)\,\into\,\CC$,
we have $|\alpha_i|=q$.
% Deligne, Weil 1, Lemme 1.7
Conjecturally, these $\alpha_i$ determine the Picard rank of $X$:

\begin{Conjecture}[Tate \cite{Tate}]
 \label{conj: tate}
   For a smooth and proper variety $X$ over $\FF_q$,
   the Picard rank $\rho(X)$ is equal to the multiplicity
   of the factor $(1-qt)$ in $P_2(t)$.
\end{Conjecture}

Although there exist more general versions of this conjecture
(see \cite{TateMotives}, for example), 
this version is sufficient for our purposes.
The following lemma is crucial for our discussion.

\begin{Lemma}
 \label{lem: supersingular}
  Let $X$ be a smooth and proper variety over $\FF_q$.
  If $X$ satisfies Conjecture \ref{conj: tate} and if
  $\Hcris{2}(\overline{X}/W)\otimes_W K$ is of slope one, then
  $\rho(\overline{X})=b_2(\overline{X})$.
\end{Lemma}

\prf
After possibly replacing $\FF_q$ by a finite extension, 
there exists a $K$-basis $\{e_i\}$ of $\Hcris{2}(X/W)\otimes_WK$ 
such that Frobenius acts as $F(e_i)=p\cdot e_i$ for all $i$.
If $q=p^r$, then $P_2(t)$ in Equation \eqref{P2} is equal to the determinant of
$({\rm id}-(F^{r})^\ast t)$ on $\Hcris{2}(X/W)\otimes_WK$,
and we conclude $P_2(t)=(1-qt)^{b_2(X)}$.
Thus, the assertion follows from Conjecture \ref{conj: tate}.
\qed\medskip

\begin{Theorem}
  \label{thm: main assuming Tate}
  If Conjecture \ref{conj: tate} holds for Enriques surfaces 
  over finite fields, then Theorem \ref{b2=rho=10} holds true.
\end{Theorem}

\prf
By Corollary \ref{cor: reduce to finite}, it suffices
to establish Theorem \ref{b2=rho=10} for 
Enriques surfaces that can be defined over finite fields.
In this special case, the claim follows from Conjecture \ref{conj: tate}
by Proposition \ref{prop: h2w is torsion} and Lemma \ref{lem: supersingular}.
\qed\medskip

\begin{Remarks} {}\quad{} 
 \label{rem: analog}
 \begin{enumerate}
 \item
 In order to establish Conjecture \ref{conj: tate} for a smooth and proper variety 
 $X$ over $\FF_q$, it suffices to 
 establish it for $X\times_{\Spec\FF_q}\Spec\FF_{q^n}$ for some $n\geq1$.
 Thus, conversely, Theorem \ref{b2=rho=10} for Enriques surfaces over $\overline{\FF}_p$ implies
 Conjecture \ref{conj: tate} for Enriques surfaces over finite fields.
 \item
 Our approach is close to the classical proof over the 
 complex numbers sketched in the introduction.
 We mention the following analogies.
 $$
  \begin{array}{l|l}
   \multicolumn{1}{c|}{\CC} & \multicolumn{1}{c}{\overline{\FF}_p} \\
   \hline
   H^2(\OO_X)=0 & H^2(W\OO_X) \mbox{ is $W$-torsion} \\
   H^{1,1}(X)=\HdR{2}(X,\CC) & \Hcris{2}(X/W)\otimes_W K \mbox{ is of slope one}\\
   \mbox{Lefschetz theorem on $(1,1)$ classes} & \mbox{Tate conjecture for divisors}
  \end{array}
 $$
 \end{enumerate}
 \end{Remarks}

\section{The Tate--conjecture for Enriques surfaces}
\label{sec: Tate}

So far, we established Theorem \ref{b2=rho=10} assuming
the Tate conjecture for divisors for Enriques surfaces over finite fields.
At the moment, it is not clear, when this conjecture
will be established in full generality, which is why we give in this section
a proof of it for Enriques surfaces to obtain an unconditional
proof of Theorem \ref{b2=rho=10}.

\subsection{The K3-like cover}
For a projective variety $X$ over a field $k$, 
we denote by $\Pic_{X/k}^\tau$
the open subgroup scheme of $\Pic_{X/k}$ that parametrizes 
divisor classes that are numerically equivalent to zero.

\begin{Theorem}[Bombieri--Mumford {\cite[Theorem 2]{bm3}}]
  \label{thm: pictau}
  If $X$ is an Enriques surface over a field $k$, then
  $\Pic^\tau_{X/k}$ is a finite group scheme of length $2$ over $k$.
\end{Theorem}

We denote by $-^D:={\cal H}om(-,\GG_m)$ Cartier duality for finite, flat,
and commutative group schemes.
Then, Theorem \ref{thm: pictau} and \cite[Proposition (6.2.1)]{Raynaud} 
(see also \cite[Section 3]{bm3} for a treatment already adapted to Enriques surfaces), 
show that, given an Enriques surface $X$ over $k$, 
the natural inclusion $\Pic^\tau_{X/k}\to\Pic_{X/k}$
gives rise to a non-trivial torsor
$$
   \pi\,:\,\widetilde{X}\,\to\,X
$$
under $(\Pic^\tau_{X/k})^D$.
In particular, $\pi$ is a finite and flat morphism of degree $2$.
Moreover, if ${\rm char}(k)\neq2$, then $\pi$ is \'etale and $\widetilde{X}$ 
is a smooth surface.
In any case, $\widetilde{X}$ is called the \emph{K3-like} cover of $X$.
The following result is a special case of \cite[Theorem 2]{Blass}, 
see also \cite[Chapter I.3]{CDL}.
 
\begin{Theorem}[Blass]
  \label{thm: K3likecover}
  If $X$ is an Enriques surface over an algebraically closed field $k$, then
  $\widetilde{X}$ is birationally equivalent to a K3 surface 
  or to $\PP^2$.
\end{Theorem}

\prf
Using the cohomological invariants in Equation \eqref{invariants} of $X$, it follows
that $\widetilde{X}$ is an integral Gorenstein surface with
$\omega_{\widetilde{X}}\iso\OO_{\widetilde{X}}$ and $\chi(\OO_{\widetilde{X}})=2$.
Let $f:Y\to\widetilde{X}$ be the minimal resolution of singularities
of the normalization of $\widetilde{X}$.

{\sc Case 1}. Assume that $\widetilde{X}$ is normal with at worst rational singularities.
Being Gorenstein, $\widetilde{X}$ has at worst rational double point singularities.
We conclude $\omega_Y\iso f^\ast\omega_{\widetilde{X}}\iso\OO_Y$
and $\chi(\OO_Y)=\chi(\OO_{\widetilde{X}})=2$ , which identifies $Y$ as a K3 surface.

{\sc Case 2}. If $\widetilde{X}$ is non-normal or normal with non-rational singularities, then
it is easy to see that $h^0(\omega_Y^{\otimes n})=0$ for all $n\geq1$.
Thus, $Y$ is of Kodaira dimension $-\infty$.
% , see \cite[Chapter I.3]{CDL} for details. 
Since $\widetilde{X}$ is not smooth,  we have ${\rm char}(k)=2$ and
$\pi$ is purely inseparable.
This implies $b_1(Y)=b_1(X)=0$ and thus, $Y$ is a rational surface, i.e., birationally
equivalent to $\PP^2$.
(We refer to \cite[Theorem 2]{Blass} for details.)
\qed\medskip

\subsection{The Tate-conjecture for Enriques surfaces over finite fields.}

\begin{Theorem}
 \label{thm: Enriques Tate}
 Enriques surfaces over finite fields satisfy Conjecture \ref{conj: tate}.
\end{Theorem} 

\prf
In order to establish Conjecture \ref{conj: tate} for a smooth and proper variety
$X$ over $\FF_q$, it suffices to establish it for $X\times_{\Spec\FF_q}\Spec\FF_{q^n}$ 
for some $n\geq1$.
By \cite[Proposition (4.3)]{TateMotives} and \cite[Theorem (5.2)]{TateMotives},
we have the following implications and equivalences:
First, if $Y\dashrightarrow X$ is a dominant and rational map between smooth and proper varieties over $\FF_q$ 
and $Y$ satisfies Conjecture \ref{conj: tate}, then so does $X$.
Second, if $Y$ and $Y'$ are a smooth, proper, and birationally equivalent varieties
over $\FF_q$, then Conjecture \ref{conj: tate} holds for $Y$ if and only 
if it holds for $Y'$.

Now, let $X$ be an Enriques surface over $\FF_q$, let $\widetilde{X}\to X$ 
be the K3-like cover, and let $Y\to \widetilde{X}$ be a resolution of singularities.
After possibly replacing $\FF_q$ by a finite extension, $Y$ is birationally equivalent 
to a K3 surface or to $\PP^2$ by Theorem \ref{thm: K3likecover}.
For $\PP^2$, Conjecture \ref{conj: tate} is trivial, and for K3 surfaces,
it is established in 
\cite{Charles}, \cite{Kim}, \cite{Madapusi}, \cite{Maulik}, \cite{Nygaard}, and \cite{NygaardOgus}.
By the above remarks and reduction steps, this implies Conjecture \ref{conj: tate} for $X$.
\qed\medskip

Combining Theorem \ref{thm: main assuming Tate} and Theorem \ref{thm: Enriques Tate},
we obtain Theorem \ref{b2=rho=10}.

\end{document}